\documentclass[10pt,draft]{amsart}
\usepackage{amsmath,amssymb,amsthm}

\begin{document}
\newtheorem{lem}{Lemma}[section]
\newtheorem{prop}{Proposition}[section]
\newtheorem{cor}{Corollary}[section]
\numberwithin{equation}{section}
\newtheorem{thm}{Theorem}[section]

\theoremstyle{remark}
\newtheorem{example}{Example}[section]
\newtheorem*{ack}{Acknowledgments}

\theoremstyle{definition}
\newtheorem{definition}{Definition}[section]

\theoremstyle{remark}
\newtheorem*{notation}{Notation}
\theoremstyle{remark}
\newtheorem{remark}{Remark}[section]

\newenvironment{Abstract}
{\begin{center}\textbf{\footnotesize{Abstract}}%
\end{center} \begin{quote}\begin{footnotesize}}
{\end{footnotesize}\end{quote}\bigskip}
\newenvironment{nome}

{\begin{center}\textbf{{}}%
\end{center} \begin{quote}\end{quote}\bigskip}

\newcommand{\triple}[1]{{|\!|\!|#1|\!|\!|}}

\newcommand{\xx}{\langle x\rangle}
\newcommand{\ep}{\varepsilon}
\newcommand{\al}{\alpha}
\newcommand{\be}{\beta}
\newcommand{\de}{\partial}
\newcommand{\la}{\lambda}
\newcommand{\La}{\Lambda}
\newcommand{\ga}{\gamma}
\newcommand{\del}{\delta}
\newcommand{\Del}{\Delta}
\newcommand{\sig}{\sigma}
\newcommand{\ome}{\omega}
\newcommand{\Ome}{\Omega}
\newcommand{\C}{{\mathbb C}}
\newcommand{\N}{{\mathbb N}}
\newcommand{\Z}{{\mathbb Z}}
\newcommand{\R}{{\mathbb R}}
\newcommand{\Rn}{{\mathbb R}^{n}}
\newcommand{\Rnu}{{\mathbb R}^{n+1}_{+}}
\newcommand{\Cn}{{\mathbb C}^{n}}
\newcommand{\spt}{\,\mathrm{supp}\,}
\newcommand{\Lin}{\mathcal{L}}
\newcommand{\SSS}{\mathcal{S}}
\newcommand{\F}{\mathcal{F}}
\newcommand{\xxi}{\langle\xi\rangle}
\newcommand{\eei}{\langle\eta\rangle}
\newcommand{\xei}{\langle\xi-\eta\rangle}
\newcommand{\yy}{\langle y\rangle}
\newcommand{\dint}{\int\!\!\int}
\newcommand{\hatp}{\widehat\psi}
\renewcommand{\Re}{\;\mathrm{Re}\;}
\renewcommand{\Im}{\;\mathrm{Im}\;}

\title{On the Existence of maximizers for a family of Restriction Theorems}

\author{Luca Fanelli}
\address{Luca Fanelli:
Universidad del Pais Vasco, Departamento de
Matem$\acute{\text{a}}$ticas, Apartado 644, 48080, Bilbao, Spain}
\email{luca.fanelli@ehu.es}

\author{Luis Vega}
\address{Luis Vega: Universidad del Pais Vasco, Departamento de
Matem$\acute{\text{a}}$ticas, Apartado 644, 48080, Bilbao, Spain}
\email{luis.vega@ehu.es}

\author{Nicola Visciglia}
\address{Nicola Visciglia: Universit\'a di Pisa, Dipartimento di Matematica,
Largo B. Pontecorvo 5, 56100 Pisa, Italy}
\email{viscigli@dm.unipi.it}


\maketitle

\begin{abstract}
We prove the existence of maximizers
for a general family of restrictions operators,
up to the end-point. We also provide some counterxamples
in the end-point case.
\end{abstract}

In the sequel we shall denote by $d\mu$ any positive measure on $\R^d_\xi$.
For every fixed $d\mu$ we define
$$T_\mu :C^0(\text{supp}(d \mu))\ni \hat h(\xi)\rightarrow \int e^{ix\cdot \xi} \hat h(\xi) d\mu
\in \mathcal C^\infty(\R^d_x)$$
Given two Banach spaces $X,Y$ we denote
by ${\mathcal L}(X,Y)$ the space of linear and continuous
operators between $X$ and $Y$.

\begin{definition}
A measure $d \mu$ on $\R^d_\xi$ satisfies the restriction condition
w.r.t.  $p\in [1,\infty]$
(shortly $(RC)_p$)
provided that
$T_\mu\in {\mathcal L}(L^2(d\mu), L^p(\R^d_x))$.
\end{definition}
\begin{definition}
Assume that $d \mu$ satisfies $(RC)_p$
then we say that there is a maximizer for
$T_\mu$ w.r.t. $p$ provided that
there exists $\hat h\in L^2(d\mu)$ such that:
$$\|\hat h\|_{L^2(d\mu)}=1$$
and
$$\|T_\mu \hat h\|_{L^p(\R^d_x)}=\|T_\mu\|_{{\mathcal L}(L^2(d\mu), L^p(\R^d_x))}.$$
\end{definition}

\begin{definition}
Assume that $ d\mu$ satisfies $(RC)_p$
then we say that $\hat h_n\in L^2(d\mu)$ is a maximizing
sequence for $T_\mu$ w.r.t. $p$ provided that:
$$\|\hat h_n\|_{L^2(d\mu)}=1$$
and
$$\lim_{n\rightarrow \infty}\|T_\mu \hat h_n\|_{L^p(\R^d_x)}=
\|T_\mu\|_{{\mathcal L}(L^2(d\mu), L^p(\R^d_x))}.$$
\end{definition}

We have the following
\begin{thm}\label{cpt}
Let $d\mu$ be a positive compactly supported measure on $\R^d_\xi$
and let $$p_0(\mu )=inf \{1\leq p \leq\infty| (RC)_p \hbox{ holds for } d\mu \}.$$
Then for every $$max \{2, p_0(\mu )\}<p \leq \infty$$
there exists a maximizer for $T_\mu$ w.r.t. $p$.
More precisely for every maximizing sequence $\hat h_n(\xi)$
for $T_\mu$ w.r.t. p,
there exists $x_n\in \R^d$ such that
$e^{ix_n\cdot \xi} \hat h_n(\xi)$ is compact in $L^2(d\mu)$.
\end{thm}

In order to treat the case $p\neq \infty$
we shall use the following general fact whose proof is inspired by
(\cite{BL},\cite{L}).

\begin{prop}\label{CC}
Let $\mathcal H$ be a Hilbert space and
$T\in {\mathcal L}({\mathcal H}, L^p(\R^d))$
for a suitable $p\in (2,\infty)$.
Let $\{h_n\}_{n\in \N}\in \mathcal H$
such that:
\begin{enumerate}
\item $\|h_n\|_{\mathcal H}=1$;
\item $\lim_{n\rightarrow \infty} \|T h_n\|_{L^p(\R^d)}=\|T\|_{{\mathcal L}({\mathcal H}, L^p(\R^d))}$;
\item
$h_n \rightharpoonup \bar h\neq 0;$
\item $T(h_n)\rightarrow T(\bar h)$ a.e. in $\R^d$.
\end{enumerate}
Then $h_n\rightarrow \bar h$ in $\mathcal H$,
in particular $\|\bar h\|_{\mathcal H}=1$ and
$\|T(\bar h)\|_{L^p(\R^d)}=\|T\|_{{\mathcal L}({\mathcal H}, L^p(\R^d))}$.
\end{prop}
\
\begin{remark}
The main difference between Proposition \ref{CC}
and Lemma 2.7 in \cite{L} is
that we only need to assume
weak convergence in the Hilbert space $\mathcal H$ for
the maximizing sequence $h_n$.
On the other hand the argument in
\cite{L} works for operators defined between general
Lebesgue spaces and not necessarily
in the Hilbert spaces framework.
\end{remark}
\begin{remark}
We shall use Proposition \ref{CC} by choosing
${\mathcal H}=L^2(d\mu)$. The main point is that in the assumptions of Proposition
\ref{CC}
we do not assume a-priori the almost everywhere convergence
of the maximizing sequence (which in our concrete context
cannot be easily checked).
\end{remark}
Next result shows that in general Theorem
\ref{cpt} cannot be extended to the end-point case $p=p_0(\mu)$.\\
For every $M>0$ we consider the compactly supported measures:\\
$$d\mu_M^1= \delta_{P_M^1},
P_M^1=\{(\xi, |\xi|^2), \xi \in \R, |\xi|\leq M\};$$
$$d\mu_M^2= \delta_{P_M^2},
P_M^2=\{(\xi, |\xi|^2), \xi \in \R^2, |\xi|\leq M\};$$
$$d\sigma_M= \frac 1{\sqrt{|\xi|}}\delta_{C_M}, C_M=\cup_{\pm}\{(\xi, \pm |\xi|), \xi \in \R^3, |\xi|\leq M\}$$
where we have denoted in general
by $\delta_S$ the flat measure on $S$.
\begin{remark}
Notice that the restriction operators associated to the measures $d\mu_M^1$,
$d\mu_M^2$, $d\sigma_M$ are strictly related
to the Strichartz estimates associated respectively to the
Schr\"odinger equation in 1-D, 2-D and to the wave equation
in 3-D (provided that the initial data
are localized in frequencies).
\end{remark}

We have the following
\begin{thm}\label{sharp}
The condition $(RC)_6$ holds
for $d\mu_M^1$ and $(RC)_{4}$ holds
for $d\mu_M^2$ and $d\sigma_M$ for every $0<M\leq \infty$.
However there are not maximizers for
$$T_{\mu_M^1}, T_{\mu_M^2}, T_{\sigma_M}$$
w.r.t. to p=6, p=4, p=4 (respectively)
provided that $M\neq \infty$.
\end{thm}

\begin{remark}
In \cite{C} it is proved
the existence of maximizers for the
restriction on the sphere ${\mathcal S}^2$
w.r.t. to $p=4$ (which turns out to be the end-point
value for the restriction on ${\mathcal S}^2$).
In the best of our knowledge this is the unique result
concerning existence of maximizers
for the end-point restriction problem on a compact manifold.
\end{remark}
\section{Proof of Theorem \ref{sharp}}

We work with $d\mu_M^1$ (the same argument works for
$d\mu_M^2$ and $d\sigma_M$).
Notice that validity of $(RC)_6$ for $d\mu_M^1$
follows from the usual Strichartz estimates
$$\|e^{it\Delta} f\|_{L^6(\R^2)}\leq C \|f\|_{L^2(\R)}.$$
Moreover the maximization problem
$$\sup_{\|\hat g\|_{L^2(d \mu_M^1)=1}}
\|T_{\mu_M^1} (\hat g)\|_{L^{6}(\R^2)}$$
is equivalent to
\begin{equation}\label{KM}
\sup_{\|h\|_{L^2(\R)=1, supp \hat h(\xi)\subset (-M, M)}} \|e^{it\Delta}h\|_{L^6(\R^2)}.
\end{equation}
On the other hand by an elementary rescaling argument
we get:
\begin{equation}\label{scaling}
\sup_{\|h\|_{L^2(\R)=1, supp \hat h(\xi)\subset (-M, M)}} \|e^{it\Delta}h\|_{L^6(\R^2)}=\sup_{\|h\|_{L^2(\R)=1}} \|e^{it\Delta}h\|_{L^6(\R^2)}.
\end{equation}
By the previous identity it is easy to deduce
that if a maximizer exists for \eqref{KM} then it is necessarily a maximizer for
\begin{equation}\label{infty}
\sup_{\|h\|_{L^2(\R)=1}} \|e^{it\Delta}h\|_{L^6(\R^2)}
\end{equation}
but this is absurd since by \cite{F}
there are no maximizers for \eqref{infty}
which are compactly supported in the Fourier variables.

\section{Proof of Proposition \ref{CC} and Theorem \ref{cpt}}

{\bf Proof of Prop \ref{CC}}
By using the Br\'ezis and Lieb Lemma (see \cite{BL}) we get:
$$\|T(h_n)-T(\bar h)\|_{L^p(\R^d)}^p= \|T(h_n)\|_{L^p(\R^d)}^p-\|T(\bar h)\|_{L^p(\R^d)}^p+o(1)$$
and by the hypothesis $(3)$ in the Proposition we get
$$\|h_n - \bar h\|_{\mathcal H}^2=
\|h_n\|_{\mathcal H}^2 - \|\bar h\|_{\mathcal H}^2 + o(1).$$
In particular since $h_n$ is by hypothesis
a maximizing sequence for $T$ we get
\begin{equation}\label{first}
\|T\|_{{\mathcal L}({\mathcal H}, L^p(\R^d))}^2
=\frac{(\|T(h_n)-T(\bar h)\|_{L^p(\R^d)}^p
+\|T(\bar h)\|_{L^p(\R^d)}^p + o(1))^\frac 2p}{\|h_n -\bar h\|^2_{{\mathcal H}}+ \|\bar h\|_{{\mathcal H}}^2
+ o(1)}\end{equation}
$$\leq \frac{(\|T(h_n)-T(\bar h)\|_{L^p(\R^d)}^2
+\|T(\bar h)\|_{L^p(\R^d)}^2 + o(1))}{\|h_n -\bar h\|^2_{{\mathcal H}}+ \|\bar h\|_{{\mathcal H}}^2 + o(1)}$$
where we have used  the inequality
$$(a + b + c)^t\leq a^t+ b^t+ c^t \hbox{ }\forall a,b,c>0$$
provided that $t\leq 1$.
The estimate above implies
\begin{equation}\label{second}\|T\|_{{\mathcal L}({\mathcal H}, L^p(\R^n))}^2
\leq \frac{( \|T\|_{{\mathcal L}({\mathcal H}, L^p(\R^n))}^2 \|h_n-\bar h\|_{L^p(\R^n)}^2
+\|T(\bar h)\|_{L^p(\R^n)}^2 + o(1))}{\|h_n -\bar h\|^2_{{\mathcal H}}+
\|\bar h\|_{{\mathcal H}}^2 + o(1)}\end{equation}
and hence
$$\|T\|_{{\mathcal L}({\mathcal H}, L^p(\R^d))}^2(
\|h_n -\bar h\|^2_{{\mathcal H}}+
\|\bar h\|_{{\mathcal H}}^2 + o(1))$$$$\leq
( \|T\|_{{\mathcal L}({\mathcal H}, L^p(\R^d))}^2 \|h_n-\bar h\|_{L^p(\R^d)}^2
+\|T(\bar h)\|_{L^p(\R^d)}^2 + o(1))$$
which is equivalent to
$$\|T\|_{{\mathcal L}({\mathcal H}, L^p(\R^d))}^2
(\|\bar h\|_{{\mathcal H}}^2 + o(1))$$$$\leq
(\|T(\bar h)\|_{L^p(\R^d)}^2 + o(1)).$$
In particular the previous estimate implies
$\|T\|_{{\mathcal L}({\mathcal H}, L^p(\R^d))}^2
\leq
\left \|T
\left (\frac{\bar h}{\|\bar h\|_{\mathcal H}}\right )\right \|_{L^p(\R^d)}^2$
and due to the definition of $\|T\|_{{\mathcal L}({\mathcal H}, L^p(\R^d))}$
it implies easily the following
\begin{equation}\label{third}
\|T\|_{{\mathcal L}({\mathcal H}, L^p(\R^d))}^2
\|\bar h\|_{{\mathcal H}}^2=
\|T(\bar h)\|_{L^p(\R^d)}^2.
\end{equation}
On the other hand by \eqref{first} we can deduce
\begin{equation}\label{second*}\|T\|_{{\mathcal L}({\mathcal H}, L^p(\R^n))}^2
\leq \frac{( \|T\|_{{\mathcal L}({\mathcal H}, L^p(\R^n))}^2 \|\bar h\|_{\mathcal H}^2
+\|T(h_n-\bar h)\|_{L^p(\R^n)}^2 + o(1))}{\|h_n -\bar h\|^2_{{\mathcal H}}+
\|\bar h\|_{{\mathcal H}}^2 + o(1)}\end{equation}
and we easily get
\begin{equation}\label{gene}
\|T\|_{{\mathcal L}({\mathcal H}, L^p(\R^d))}^2
(\|h_n-\bar h\|_{{\mathcal H}}^2 + o(1))$$$$\leq
(\|T(h_n-\bar h)\|_{L^p(\R^d)}^2 + o(1)).\end{equation}
Notice that
either $\|h_n-\bar h\|_{\mathcal H}=o(1)$ (and in this case we can conclude)
or (up to subsequence)
$$\inf_{n\in \N} \|h_n-\bar h\|_{\mathcal H}\geq \epsilon_0>0.$$
In particular by \eqref{gene} we get
$$\|T\|_{{\mathcal L}({\mathcal H}, L^p(\R^d))}^2
\leq
\left\|T\left (\frac{h_n-\bar h}
{\|h_n-\bar h\|_{{\mathcal H}}}\right )\right \|_{L^p(\R^d)}^2 + o(1))$$
which by definition of
$\|T\|_{{\mathcal L}({\mathcal H}, L^p(\R^d))}$
necessarily implies
$$\|T\|_{{\mathcal L}({\mathcal H}, L^p(\R^d))}^2
=
\left\|T\left (\frac{h_n-\bar h}
{\|h_n-\bar h\|_{{\mathcal H}}}\right )\right \|_{L^p(\R^d)}^2 + o(1))$$
and equivalently
\begin{equation}\label{fourth}
\|T\|_{{\mathcal L}({\mathcal H}, L^p(\R^d))}^2
\|h_n-\bar h\|_{{\mathcal H}}^2=
\|T(h_n-\bar h)\|_{L^p(\R^d)}^2 + o(1)).
\end{equation}
By combining the first identity in \eqref{first}
with \eqref{third} and \eqref{fourth} we get
\begin{equation}\label{first*}
1
=\frac{(\|h_n-\bar h\|_{\mathcal H}^p
+\|\bar h\|_{\mathcal H}^p + o(1))^\frac 2p}{\|h_n -\bar h\|^2_{{\mathcal H}}+ \|\bar h\|_{{\mathcal H}}^2
+ o(1)}\end{equation}
i.e.
$$(\|h_n-\bar h\|_{\mathcal H}^p
+\|\bar h\|_{\mathcal H}^p)^\frac 2p=
\|h_n -\bar h\|^2_{{\mathcal H}}+ \|\bar h\|_{{\mathcal H}}^2 +o(1).$$
Since we are assuming $p\in (2, \infty)$ it is easy to deduce
by a convexity argument that the previous inequality
implies
$\|\bar h\|_{{\mathcal H}}=1$ and $\|h_n -\bar h\|^2_{{\mathcal H}}=o(1)$
(actually we have excluded the possibility
$\|\bar h\|_{{\mathcal H}}=0$ and $\|h_n -\bar h\|^2_{{\mathcal H}}=1+o(1)$
since by assumption $\bar h\neq 0$).
\hfill$\Box$

{\bf Proof of Thm \ref{cpt}}
\\
{\bf The case $p\neq \infty$}\\
\\
Let $\hat h_n\in L^2(d\mu)$ be a maximizing sequence
for $T_\mu$ w.r.t. $p$ (where $p$ is as in the assumptions).\\
{\em {\bf First step}:
there is a sequence $x_n\in \R^d$  such that
$\hat g_n(\xi)=e^{ix_n \cdot \xi} \hat h_n(\xi)$
has a weak limit different from zero in $L^2(d\mu)$}
\\
\\
In order to verify this property
we prove that there is $x_n$ such that
$$T_\mu ((e^{ix_n\cdot \xi} \hat h_n(\xi))=\tau_{x_n} T_\mu (\hat h_n(\xi))$$
has a weak limit different from zero
(here $\tau_y$ denotes the translation of vector $y$).
Notice that
by definition we have
\begin{equation}\label{1}
\|T_\mu \hat h_n\|_{L^p(\R^n_x)}\rightarrow \|T_\mu\|_{{\mathcal L}(L^2(d\mu), L^p(\R^n_x))}>0.
\end{equation}
By using the $(RC)_{\bar p}$ condition for
a suitable
$p_0(\mu)< \bar p <p$ we get
$$\|T_\mu \hat h_n \|_{L^{\bar p}(\R^n_x)}
\leq \|T\|_{{\mathcal L}(L^2(d\mu), L^{\bar p}(\R^n_x))} \|\hat h_n(\xi)\|_{L^2(d\mu)}$$
and hence
\begin{equation}\label{2}
\sup_{n\in \N} \|T_\mu \hat h_n \|_{L^{\bar p}(\R^d_x)}\equiv S<\infty.
\end{equation}
Next notice that we have the following inequality:
$$\|T_\mu \hat h_n\|_{L^p(\R^d_x)}\leq \|T_\mu \hat h_n\|_{L^{\bar p}(\R^d_x)}^\theta
\|T_\mu \hat h_n\|_{L^\infty(\R^d_x)}^{1-\theta}$$
where
$\frac 1p=\frac{\theta}{\bar p}$.
By combining this fact with \eqref{1} and \eqref{2} we deduce
\begin{equation}\label{fine}
\|T_\mu \hat h_n\|_{L^\infty(\R^d_x)}\geq \epsilon_0>0.
\end{equation}
Notice also that we have (by compactness of the support of $d\mu$)
$$\|T_\mu \hat h_n\|_{L^\infty(\R^d_x)}
\leq \|\hat h\|_{L^2(d\mu)} \sqrt{\|d\mu\|}$$
and
$$\|\nabla_x T_\mu \hat h_n\|_{L^\infty(\R^d_x)}
= \|T_\mu (i\xi \hat h_n)\|_{L^\infty(\R^d_x)}$$$$\leq \|\xi \hat h_n\|_{L^2} \sqrt{\|d\mu\|}
\leq \sqrt {\|d\mu\|} \left (\sup_{\xi\in supp (\mu)} |\xi|\right )\|\hat h_n\|_{L^2(d\mu)} $$
(where $\|d\mu\|=\int 1 d\mu$).
Hence
\begin{equation}\label{unifconv}
\sup_{n\in \N} \|T_\mu \hat h_n\|_{W^{1,\infty}(\R^d_x)}<\infty.
\end{equation}
By \eqref{fine} there exist $x_n$ such that
\begin{equation*}
|T_\mu \hat h_n(x_n)|\geq \epsilon_0>0\end{equation*}
and hence
\begin{equation}\label{tau}|\tau_{x_n} T_\mu \hat h_n(0)|\geq \epsilon_0>0.\end{equation}
On the other hand by \eqref{unifconv} we get
\begin{equation*}
\|\tau_{x_n} T_\mu \hat h_n\|_{W^{1,\infty}(B(0,1))}\end{equation*}
are uniformly bounded and hence by the Ascoli-Arzel\'a Theorem
$$\tau_{x_n} (T_\mu \hat h_n(\xi))$$ has an uniform limit in $B(0,1)$.
By \eqref{tau} the limit
has to be different from zero.
\\
{\em{\bf Second step}: conclusion of the proof}
\\
\\
Notice that
$$\|\hat g_n(\xi)\|_{L^2}=1$$
and
$$\|T_\mu (\hat g_n)\|_{L^p(\R^n_x)}=
\|T_\mu (\hat h_n)\|_{L^p(\R^n_x)}.$$
Hence $\hat {g_n}$ is a maximizing sequence for $T_\mu$.
On the other hand
by the previous step
it is easy to check that all the hypothesis
of Proposition \ref{CC} are satisfied if we choose
$T=T_\mu$, ${\mathcal H}=L^2(d\mu)$ and we fix as a maximizing sequence
$\hat g_n$.
\\
{\bf  The case $p=\infty$}
\\
\\
Following the computations done above we have that
$$
\lim_{n\rightarrow \infty}
\|T_\mu \hat h_n\|_{L^\infty (\R^n_x)}=\|T_\mu\|_{{\mathcal L}(L^2(d\mu), L^\infty(\R^d_x))}
$$
and moreover
\begin{equation}\label{last}
\sup_{n\in \N}\|T_\mu \hat h_n\|_{W^{1,\infty}(\R^d_x)}<\infty.
\end{equation}
In particular there is a sequence
$x_n\in \R^d$ such that
$$\lim_{n\rightarrow \infty}
\|T_\mu \hat h_n(x_n)\|=\|T_\mu\|_{{\mathcal L}(L^2(d\mu), L^\infty(\R^d_x))}.$$
As in the previous case we introduce
$\hat g_n=e^{ix_n\cdot \xi}\hat h(\xi)$
and it is easy to deduce that
$\hat g_n$ is still  maximizing sequence
with the extra property that
\begin{equation}\label{anal}
\lim_{n\rightarrow \infty}
|T_\mu (\hat g_n)(0)|=\|T_\mu\|_{{\mathcal L}(L^2(d\mu), L^\infty(\R^d_x))}.
\end{equation}
By the Ascoli-Arzel\'a theorem (that can be applied due to \eqref{last})
in conjunction with \eqref{anal} we
conclude that if $\bar g$ is the the weak limit
of $\hat g_n$ in $L^2(d\mu)$ then necessarily
$$\|T_\mu (\bar g)\|_{L^\infty(\R^n_x)}
\geq |T_\mu \bar g(0)|=\|T_\mu\|_{{\mathcal L}(L^2(d\mu), L^\infty(\R^d_x))}.$$
On the other hand by semicontinuity of the norm
$L^2(d\mu)$ we have that
$\|\bar g\|_{L^2(d\mu)}\leq 1$. By combining this fact with the definition
of
$\|T_\mu\|_{{\mathcal L}(L^2(d\mu), L^\infty(\R^d_x))}$
we easily deduce that $\|\bar g\|_{L^2(d\mu)}=1$
and hence $\hat g_n$ is compact in $L^2(d\mu)$.

\end{document}